\input amstex
\documentstyle{amsppt}
\mag1095
\hfuzz5.4pt

\hyphenation{Berg-man Sobo-lev Toep-litz}

\define\CC{{\Bbb C}}
\define\RR{{\Bbb R}}
\define\BB{{\Bbb B}}
\define\KK{\bold K}
\define\spr#1{\langle#1\rangle}
\define\emph#1{{\sl #1\/}}
\define\hol{_{\text{hol}}}
\define\harm{_{\text{harm}}}
\define\oOm{{\overline\Omega}}
\define\pOm{{\partial\Omega}}
\define\ada{A^2_\alpha}
\define\adah{A^2_{\alpha*}}
\define\dbar{\overline\partial}
\define\psdo/{$\Psi$DO}
\define\hdva{H^2(\pOm)}
\redefine\Re{\operatorname{Re}}
\redefine\Im{\operatorname{Im}}
\define\ord{\operatorname{ord}}
\define\Ker{\operatorname{Ker}}
\define\Ran{\operatorname{Ran}}
\loadbold
\define\jedna{\boldkey1}
\define\diag{\operatorname{diag}}

\define\ccz{\CC[z_1,\dots,z_d]}
\define\MM{{\Cal M}}
\define\MMp{{\MM^\perp}}
\define\HH{\Cal H}
\define\SS{\Cal S}
\define\adac{A^2_{\alpha\circ}}
\define\gto/{generalized Toeplitz operator}
\define\sprb#1{\spr{#1}_\pOm}
\define\spro#1{\spr{#1}_\Omega}
\define\cD{\Cal D}
\define\dR{\dot{\RR}{}}
\define\Tb{T^\bullet}
\redefine\ell{_{\text{ell}}}
\define\rdv{R_{\partial V}}
\define\rdvpb{\rdv\Pi}
\define\sibo{\Sigma_{\pOm|V}}
\define\sprdb#1{\spr{#1}_{\partial\BB}}

\define\akov{A^2_{\alpha+k,*}(\Omega\cap V)}
\define\pOcv{{\pOm\cap V}}
\define\bcv{{\BB^d\cap V}}
\define\pbcv{{\partial\BB^d\cap V}}
\define\akbv{A^2_{\alpha+k,*}(\bcv)}
\define\cV{{\Cal V}}
\define\cL{{\Cal L}}
\define\cB{{\Cal B}}
\define\ako{A^2_{\alpha+k,*}(\Omega)}
\define\cN{\Cal N}
\define\cc{\Cal C}

\define\cT{\Cal T}
\define\cH{\Cal H}

\let\oldtag\tag
\newcount\tagno \tagno=0
\def\tag#1{\ifx#1\undefined\relax\global\advance
        \tagno by1\global\edef#1{\the\tagno}
        \else\errmessage{TAG \string#1\space is already defined}\fi \oldtag#1}
\def\tagdef#1{\ifx#1\undefined\relax\global\advance
        \tagno by1\global\edef#1{\the\tagno}
        \else\errmessage{TAG \string#1\space is already defined}\fi}
\newcount\prono \prono=0
\def\prodef#1{\ifx#1\undefined\relax\global\advance
        \prono by1\global\edef#1{\the\prono}\else
        \errmessage{TAG \string#1\space is already defined}\fi #1}
\define\chk#1#2{\ifx#1#2\relax\else
                \errmessage{Change \string#1\space to #2!}\fi}
\define\prPY{4}
\define\prtFF{18}

\newcount\refno \refno=0
\define\refdef#1{\ifx#1\undefined\relax\global\advance
        \refno by1\global\edef#1{\the\refno}
        \else\errmessage{REF \string#1\space is already defined}\fi}
\define\refd#1#2{\ifkonec\item{[#1]} #2.\par\else\refdef#1\fi}
\newif\ifkonec \konecfalse

\refd\ArvP {W. Arveson: {\it $p$-summable commutators in dimension~$d$,\/}
 J.~Operator Theory {\bf 54} (2005), 101--117}

\refd\BdMA {L. Boutet de Monvel: {\it Boundary problems for pseudo-differential
operators,\/} Acta Math. {\bf 126} (1971), 11--51}

\refd\BdMi {L. Boutet de Monvel: {\it On the index of Toeplitz operators 
in several complex variables,\/} Invent. Math. {\bf 50} (1979), 249--272}

\refd\BdMG {L. Boutet de Monvel, V. Guillemin: {\it The spectral theory of
Toeplitz operators,\/} Ann. Math. Studies, vol.~99, Princeton University 
Press, Princeton, 1981}

\refd\BdMS {L. Boutet de Monvel, J. Sj\"ostrand: {\it Sur la singularit\'e des
noyaux de Bergman et de Szeg\"o,\/} Ast\'erisque {\bf 34--35} (1976), 123--164}
 
\refd\BdMT {L. Boutet de Monvel: {\it Symplectic cones and Toeplitz
operators,\/} Multidimensional Complex Analysis and Partial Differential
Equations, pp.~15--24, Contemporary Math., vol.~205, Amer. Math. Soc.,
Providence, 1997}

\refd\Bea {F. Beatrous: {\it Estimates for derivatives of holomorphic
functions in pseudoconvex domains,\/} Math. Z. {\bf 191} (1986), 91--116}

\refd\Doug {R. Douglas: {\it A~new kind of index theorem,\/} Analysis,
geometry and topology of elliptic operators, pp.~369-382, World Sci. Publ.,
Hackensack, 2006}

\refd\DW {R.G. Douglas, K.Wang: {\it Essential normality of cyclic submodule
generated by any polynomial,\/} preprint, 2011, arXiv:1101.0774}

\refd\Ejfa {M. Engli\v s: {\it Toeplitz operators and weighted Bergman
kernels,\/} J.~Funct. Anal. {\bf 255} (2008), 1419--1457}

\refd\Epay {M. Engli\v s: {\it Weighted Bergman kernels for logarithmic
weights,\/} Pure Appl. Math. Quarterly (Kohn special issue) {\bf6} (2010),
781--813} 

\refd\E {M. Engli\v s: {\it Analytic continuation of weighted Bergman
kernels,\/} J.~Math. Pures Appl. {\bf 94} (2010), 622--650}

\refd\EZ {M. Engli\v s, G. Zhang: {\it Hankel operators and the Dixmier
trace on strictly pseudoconvex domains,\/} Docum. Math. {\bf 15} (2010),
601--622}

\refd\FX {Q. Fang, J. Xia: {\it Essential normality of polynomial-generated
submodules: Hardy space and beyond,\/} preprint, 2011, available~at 
http://www.acsu.buffalo.edu/$^\sim$jxia/Preprints/poly.pdf}

\refd\GrauRem {H. Grauert, R.Remmert, {\it Coherent analytic sheaves,\/}
Springer, Berlin, 1984}

\refd\GrSj {A. Grigis, J. Sj\"ostrand, {\it Microlocal analysis for
diferential operators,\/} London Math. Soc. Lecture Notes, vol.~196,
Cambridge Univ. Press, Cambridge, 1994}

\refd\GW {K. Guo and K. Wang, {\it Essentially normal Hilbert modules and
K-homology,\/} Math. Ann. {\bf 340} (2008), 907--934}

\refd\Hcpam {L. H\"ormander: {\it Pseudo-differential operators,\/} Comm.
Pure Appl. Math. {\bf18} (1965), 501--517}

\refd\Horm {L. H\"ormander, {\it The analysis of linear partial differential
operators, vol.~I--IV,\/} Grund\-leh\-ren der mathematischen Wissenschaften,
Springer-Verlag, 1985}

\refd\KS {M. Kennedy, O.M. Shalit: {\it Essential normality and the
decomposability of algebraic varieties,\/} New~York J. Math. {\bf 18}
(2012), 877--890}

\refd\LM {J.-L. Lions, E. Magenes, {\it Probl\`emes aux limites non
homog\`enes et applications,\/} vol.~1, Dunod, Paris, 1968}

\refd\MeSj {A. Melin, J. Sj\"ostrand: {\it Fourier integral operators
with complex valued phase functions,\/} Fourier integral operators and
partial differential equations, Lecture Notes Math., vol.~459, pp.~120-223,
Springer Verlag, Berlin-Heidelberg, 1975}

\refd\Sh {O.M. Shalit: {\it Operator theory and function theory in
Drury-Arveson space and its quotients,\/} preprint, 2013, arxiv:1308.1081}

\refd\Strich {R. Strichartz: {\it A~functional calculus for elliptic
pseudo-differential operators,\/} Amer. J. Math. {\bf 94} (1972), 711--722}

\refd\Trev {F. Treves, {\it Introduction to pseudodifferential and Fourier
integral operators,\/} Plenum, New~York, 1980}


\topmatter
\title Geometric Arveson-Douglas conjecture\endtitle
\rightheadtext{Arveson-Douglas conjecture}
\leftheadtext{M.~Engli\v s, J.~Eschmeier}
\author Miroslav Engli\v s, J\"org Eschmeier\endauthor
\address Mathematics Institute, \v Zitn\' a 25, 11567~Prague~1, Czech Republic
 {\rm and} Mathematics Institute, Silesian University at Opava,  
 Na~Rybn\'\i\v cku~1, 74601~Opava, Czech Republic\endaddress 
\email englis{\@}math.cas.cz\endemail
\address Fachrichtung Mathematik, Universit\"at des Saarlandes,
 Postfach~15~11~50, D-66041~Saarbr\"ucken, Germany\endaddress
\email eschmei{\@}math.uni-sb.de\endemail
\subjclass Primary 47A13; Secondary 32W25, 47B35\endsubjclass
\keywords Arveson-Douglas conjecture, generalized Toeplitz operator\endkeywords
\abstract
We~prove the Arveson-Douglas essential normality conjecture for graded
Hilbert submodules that consist of functions vanishing on a given homogeneous
subvariety of the ball, smooth away from the origin. Our~main tool is the
theory of \gto/s of Boutet de Monvel and Guillemin. \endabstract
\endtopmatter

\document

\head 1. Introduction\endhead
Let $\BB^d$ be the unit ball in $\CC^d$, $d\ge1$. The~Drury-Arveson space
$H^2_d$ consists of all holomorphic functions
$f(z)=\sum_\nu f_\nu z^\nu$ on $\BB^d$ such that
$$ \|f\|^2_{DA} := \sum_\nu |f_\nu|^2 \frac{\nu!}{|\nu|!} < \infty ,   $$
equipped with the corresponding norm and inner product. The~operators
$$ M_{z_j}: f(z) \mapsto z_j f(z)   $$
of~multiplication by the coordinate functions are bounded on~$H^2_d$, and
commute with each other. This endows $H^2_d$ with the structure of a module
over the polynomial ring $\ccz$, a~polynomial $p$ corresponding to the
operator $M_p=p(M_{z_1},\dots,M_{z_d})$ of multiplication by $p$ on~$H^2_d$.
If~$\MM\subset H^2_d$ is a (closed) subspace invariant under all~$M_{z_j}$,
$j=1,\dots,d$, we~can therefore consider the restrictions $M_{z_j}|_\MM$,
which are commuting bounded linear operators on~$\MM$, as~well as the
compressions 
$$ S_j := P_\MMp M_{z_j} |_\MMp, \qquad j=1,\dots,d,   $$
of~the $M_{z_j}$ to the orthogonal complement $\MMp=H^2_d\ominus\MM$, which are
commuting bounded linear operators on~$\MMp$.

The~following conjecture was originally made by Arveson~\cite{\ArvP} with $d$
in the place of $\dim Z(p)$, and refined to the current form by Douglas~\cite
{\Doug}.~\footnote{In~both cases, it~was also formulated for the more general
case of modules $\MM$ in $H^2_d\otimes\CC^N$ generated by $\CC^N$-valued
homogeneous polynomials, with some finite~$N\ge1$.}

\proclaim{Arveson-Douglas Conjecture} Assume $\MM$ is generated, as~a~module,
by~finitely many homogeneous polynomials $p_1,\dots,p_m\in\ccz$.
Then the commutators $[S_j,S^*_k]$, $j,k=1,\dots,d$, belong to the Schatten
class $\SS^q$ for all $q>\dim Z(p)$, where $\dim Z(p)$ is the complex dimension
of the zero-set $Z(p)\equiv Z(p_1,\dots,p_m)$ of the polynomials
$p_1,\dots,p_m$. \endproclaim 

The~Arveson conjecture, and in some cases also its refined version due
to~Douglas, have so far been proved in various special settings: by~Arveson
himself \cite{\ArvP} when $p_1,\dots,p_m$ are monomials; by~Guo and~Wang
\cite{\GW} for $m=1$ or $d\le3$; by~Douglas and~Wang \cite{\DW} when $m=1$
and $\MM$ is a submodule of the Bergman space $L^2\hol(\BB^d)$ on~$\BB^d$
(instead~of~$H^2_d$) generated by an arbitrary, not necessarily homogeneous
polynomial~$p$; by~Fang and Xia~\cite{\FX} for submodules of the same type
in certain weighted (Sobolev-)Bergman spaces on~$\BB^d$, which included
$L^2\hol(\BB^d)$ as~well as the Hardy space $H^2(\partial\BB^d)$ on~$\BB^d$,
but not~$H^2_d$ (unless $d=1$); by~Kennedy and~Shalit \cite{\KS} when
$p_1,\dots,p_m$ are homogeneous polynomials such that the linear spans of
$Z(p_1),\dots,Z(p_m)$ in~$\CC^d$ have mutually trivial intersections; etc.
See~the recent survey paper by Shalit \cite{\Sh} for some more details and
further information, as~well as the original paper by Douglas~\cite{\Doug}
for more on the motivation and applications to $K$-homology and index theory.

There is also a reformulation of (a~weaker version~of) the~Arveson-Douglas
conjecture in terms of varieties. Namely, denote by $I(p)$ the ideal in $\ccz$
generated by $p_1,\dots,p_m$; then $\MM$ is the closure of $I(p)$ in~$H^2_d$,
and $I(p)$ is a \emph{homogeneous} (or~\emph{graded})~ideal, meaning that 
$$ I(p)=\sum\nolimits_{k\ge0}^\oplus(I(p)
 \cap(\text{homogeneous polynomials of degree }k)).  $$
Denoting for any ideal $J$ in~$\ccz$~by
$$ Z(J) := \{z\in\CC^d: q(z)=0 \;\forall q\in J \}   $$
the~zero set of~$J$, we~then have $Z(p)=Z(I(p))$, which is a homogeneous
variety in~$\CC^d$, i.e.~$z\in Z(p)$, $t\in\CC$ implies $tz\in Z(p)$.
Conversely, for~any subset $X\subset\CC^d$,
$$ I(X) := \{q\in\ccz: q(z)=0 \;\forall z\in X \}   $$
is~an ideal in~$\ccz$, which is homogeneous if $X$~is. The~correspondences
$J\mapsto Z(J)$, $X\mapsto I(X)$ are not one-to-one: one~always has
$I(Z(J))\supset J$, with equality if and only if $J$ is a radical ideal,
i.e.~$J=\sqrt J$ where $\sqrt J:=\{q\in\ccz: q^n\in J \text{ for some }
n=1,2,\dots\}$; also, $Z(J_1)=Z(J_2)$ if~and only~if $\sqrt{J_1}=\sqrt{J_2}$
(this~is Hilbert's Nullstellensatz). Specializing to modules generated by
radical ideals, we~thus get the following ``geometric version'' of the
Arveson-Douglas conjecture~\cite{\KS}.

\proclaim{Geometric Arveson-Douglas conjecture}
Let $V$ be a homogeneous variety in~$\CC^d$ and $\MM=\{f\in H^2_d: f(z)=0
\text{ for all }z\in V\cap\BB^d\}$. Then $[S_j,S^*_k]\in\SS^q$ for all
$q>\dim_\CC V$.    \endproclaim

As~already mentioned in passing, one~can consider the above conjectures not
only for~$H^2_d$, but also for other spaces of holomorphic functions on~$\BB^d$
on which the multiplication operators $M_{z_j}$, $j=1,\dots,d$, act boundedly.
These include the (weighted Bergman) spaces
$$ \ada(\BB^d) \equiv \ada := L^2\hol(\BB^d,d\mu_\alpha)   $$
of~holomorphic functions on~$\BB^d$ square-integrable with respect to the
probability measure
$$ d\mu_\alpha(z) := \frac{\Gamma(\alpha+d+1)}{\Gamma(\alpha+1)\pi^d}
 (1-|z|^2)^\alpha \, dz,  \qquad \alpha>-1,   $$
where $dz$ denotes the Lebesgue volume on $\CC^d$ and the restriction
on $\alpha$ ensures that these spaces are nontrivial (and~contain all
polynomials). In~terms of the Taylor coefficients $f(z)=\sum_\nu f_\nu z^\nu$,
the~norm in $\ada$ is given~by 
$$ \|f\|_\alpha^2 = \sum_\nu |f_\nu|^2 \frac{\nu!\,\Gamma(d+\alpha+1)}
 {\Gamma(|\nu|+d+\alpha+1)}.   \tag\tTF   $$
The~right-hand side makes actually sense and is positive-definite for all
$\alpha>-d-1$, and we can thus extend the definition of $\ada$ also to $\alpha$
in this range; in~particular, this will give, in~addition to the weighted
Bergman spaces for $\alpha>-1$ (including the ordinary --- i.e.~unweighted ---
Bergman space $L^2\hol(\BB^d)$ for $\alpha=0$), also the Hardy space
$$ A^2_{-1} = H^2(\partial\BB^d,d\sigma)   $$
with respect to the normalized surface measure $d\sigma$ on~$\partial\BB^d$
for $\alpha=-1$, as~well as the Drury-Arveson space 
$$ A^2_{-d} = H^2_d   $$
for $\alpha=-d$. Furthermore, passing from (\tTF) to the equivalent norm
$$ \|f\|^2_{\alpha\circ} := \sum_\nu \frac{|f_\nu|^2}{(|\nu|+1)^{d+\alpha}}
 \, \frac{\nu!}{|\nu|!} ,  \tag\tAC  $$
one~can even define the corresponding spaces $\adac$ for any
real~$\alpha$, with $\adac=\ada$ (as~sets, with equivalent norms)
for~$\alpha>-d-1$ (hence, in~particular, $A^2_{-d,\circ}=H^2_d$ for $\alpha=-d$,
$A^2_{-1,\circ}=H^2(\partial\BB^d)$ for $\alpha=-1$, and $\adac=\ada$
for $\alpha>-1$). Actually, $\adac$~are precisely the subspaces
of holomorphic functions 
$$ \adac = W^{-\alpha/2}\hol(\BB^d) := \{f\in W^{-\alpha/2}(\BB^d): \;
 f\text{ is holomorphic on }\BB^d \}   $$
in~the Sobolev spaces $W^{-\alpha/2}(\BB^d)$ on $\BB^d$ of order
$-\frac\alpha2$, for any real $\alpha$. The~coordinate
multiplications~$M_{z_j}$, $j=1,\dots,d$, are~continuous on $\adac$ for
any~$\alpha\in\RR$, and one can consider the Arveson-Douglas conjecture
in this setting.

Our~main result is the proof of the geometric variant of the Arveson-Douglas
conjecture --- that~is, proof of the Arveson-Douglas conjecture for
subspaces $\MM$ generated by a radical homogeneous ideal --- in~all these
settings for smooth submanifolds. 

\proclaim{Main Theorem} Let $V$ be a homogeneous variety in $\CC^d$ such that
$V\setminus\{0\}$ is a complex submanifold of $\CC^d\setminus\{0\}$ of
dimension~$n$, $\alpha\in\RR$, and $\MM$ the subspace in~$\adac$, or~in $\ada$
if~$\alpha>-d-1$, of functions vanishing on~$V\cap\BB^d$.
Then $[S_j,S^*_k]\in\SS^q$, $j,k=1,\dots,d$, for all $q>n$.    \endproclaim

Our~method of proof relies on two ingredients: the~results of Beatrous about
restrictions of functions in $\adac$ to submanifolds~\cite{\Bea}, and the
theory of Boutet de Monvel and Guillemin of Toeplitz operators on the Hardy
space with pseudodifferential symbols (so-called ``\gto/s'')
\cite{\BdMG}~\cite{\BdMi}. It~actually turns out that the Boutet de Monvel
and Guillemin theory can also be used to replace the results of Beatrous
from~\cite{\Bea}, at~least those that we need here. The~required prerequisites
about the \gto/s of Boutet de Monvel and Guillemin are reviewed in Section~2,
and those about restrictions to submanifolds in Section~3.
With~these tools it is possible to prove a variant of our main theorem with $V$
a (not~necessarily homogeneous) complex submanifold of~$\BB^d$ intersecting
$\partial\BB^d$ transversally; we~do this in Section~4. The~proof of Main
Theorem, which builds on the same ideas but with some additional
technicalities, is~given in Section~5.

\head 2. Generalized Toeplitz operators\endhead
Let~$\Omega$ be a bounded strictly pseudoconvex domain with smooth
(i.e.~$C^\infty$) boundary in~$W$, where $W$ is either $\CC^n$ or,
more generally, a~complex manifold of dimension $n$; an~example is $W$
a~complex submanifold of dimension $n$ in $\CC^d$, $d>n$, 
and $\Omega=W\cap\BB^d$.
(One~could even allow $W$ to be a complex analytic variety of dimension~$n$
with singularities in $\Omega$ but not on~$\pOm$, an~example being
a~homogeneous complex cone of dimension $n$ in~$\CC^d$, $d>n$, again with
$\Omega=W\cap\BB^d$; see~\S2i in~\cite{\BdMi}.)
We~fix a positively-signed ``defining function'' $\rho$ for~$\Omega$, 
i.e.~a~function smooth on the closure $\oOm$ of $\Omega$ such that
$\rho>0$ on $\Omega$ and $\rho=0$, $\nabla\rho\neq0$ on~$\pOm$;
in~the example above, we~can take $\rho(z)=1-|z|^2$.

Let~$L^2(\pOm)$ be the Lebesgue space on the boundary $\pOm$ with respect
to the surface measure (i.e.~the $(2n-1)$-dimensional Hausdorff measure)
$d\lambda$~on~$W$; we~will denote the inner product and norm in
$L^2(\pOm)$ by~$\sprb{\cdot,\cdot}$ and~\hbox{$\|\cdot\|_\pOm$}, respectively,
and similarly by $\spro{\cdot,\cdot}$ the inner product in~$L^2(\Omega,dz)$.
\footnote{In~the case of $W$ a manifold, we~use the surface measure and
volume with respect to some chosen Riemannian metric on~$W$.}
Let~$\dbar$ denote the usual Cauchy-Riemann operator on~$W$, and $\dbar{}^*$
its (formal) adjoint with respect to some fixed smooth Hermitian metric on~$W$;
the~harmonic functions on~$W$ are then those annihilated by the Laplacian
$\Delta:=-\dbar{}^*\dbar$. 
The~Hardy space $\hdva$ is the subspace in $L^2(\pOm)$ of functions
whose Poisson extension into $\Omega$ is not only harmonic but holomorphic;
or,~equivalently, the~closure in $L^2(\pOm)$ of $C^\infty\hol(\pOm)$,
the~space of boundary values of all the functions in $C^\infty(\oOm)$ that
are holomorphic on~$\Omega$. We~will also denote by $W^s(\pOm)$, $s\in\RR$,
the~Sobolev spaces on~$\pOm$, and by $W^s\hol(\pOm)$ the closure of
$C^\infty\hol(\pOm)$ in~$W^s(\pOm)$. The~Poisson extension operator
$$ \KK: C^\infty(\pOm)\to C^\infty(\oOm),
 \qquad \Delta\KK u=0 \text{ on }\Omega, \qquad \KK u|_\pOm = u,   $$
then extends to a bounded operator from $W^s(\pOm)$ onto $W^{s+1/2}\harm
(\Omega)$, the~subspace of harmonic functions in the Sobolev space $W^{s+1/2}
(\Omega)$ on~$\Omega$, and~from $W^s\hol(\pOm)$ onto~$W^{s+1/2}\hol(\Omega)$,
the~subspace of holomorphic functions in~$W^{s+1/2}(\Omega)$. The~operator of
taking the boundary values (or~``trace'')
$$ \gamma: C^\infty(\oOm)\to C^\infty(\pOm), \qquad \gamma f:=f|_\pOm,  $$
which acts from $W^s(\Omega)$ onto $W^{s-1/2}(\pOm)$ for $s>\frac12$
(this is the Sobolev trace theorem), similarly extends to a bounded map from
$W^s\harm(\Omega)$ onto~$W^{s-1/2}(\pOm)$ and from $W^s\hol(\Omega)$
onto~$W^{s-1/2}\hol(\pOm)$, for~any $s\in\RR$, which is the right inverse
to~$\KK$. On~harmonic and holomorphic functions, $\gamma$~and $\KK$ are thus
mutual inverses, establishing isomorphisms $W^s\harm(\Omega)\leftrightarrow
W^{s-1/2}(\pOm)$ and $W^s\hol(\Omega)\leftrightarrow W^{s-1/2}\hol(\pOm)$,
for~any real~$s$. See~e.g.~Lions and Magenes~\cite{\LM}, Chapter~2, Section~7.3
for the proofs and further details.
\footnote{Or~page~29, fifth paragraph, in~\cite{\BdMA}; as~mentioned explicitly
on p.~13 there, this applies also to the case of manifolds~$\Omega$, not~only
to domains in~$\RR^n$.}

As~usual, by~a~classical pseudodifferential operator (or~\psdo/ for short)
on~$\pOm$ of order~$m$ we~will mean a pseudodifferential operator whose total
symbol in any local coordinate system has an asymptotic expansion
$$ p(x,\xi) \sim \sum_{j=0}^\infty p_{m-j}(x,\xi)   $$
where $p_{m-j}$ is $C^\infty$ in $x,\xi$ and positive homogeneous of degree
$m-j$ in $\xi$ for $|\xi|>1$. Here $m$ can be any real number, and the symbol
``$\sim$'' means that the difference between $p$ and $\sum_{j=0}^{k-1}p_{m-j}$ 
should belong to the H\"ormander class $S^{m-k}$, for each $k=0,1,2,\dots$;
see~\cite{\Horm}. The~space of all such operators will be denoted~$\Psi^m$.
An~operator in $\Psi^m$ maps $W^s(\pOm)$ into $W^{s-m}(\pOm)$ for any
$s\in\RR$. Unless explicitly stated otherwise, all \psdo/s in this paper
will be classical. 

If~$A\in\Psi^m$, $m<0$, is~elliptic, i.e.~its principal symbol $\sigma(A)
(x,\xi)=a_m(x,\xi)$ does not vanish for $\xi\neq0$, and is positive selfadjoint
as an operator on $L^2(\pOm)$ (i.e.~$\spr{Au,u}>0$ for all $u\in L^2(\pOm)$,
$u\neq0$), then $A$ is compact and its spectrum consists of isolated
eigenvalues $\lambda_1>\lambda_2>\dots>0$ of finite multiplicity, 
so~one can define the power $A^z$ for any $z\in\CC$ by the spectral theorem.
Similarly for positive (i.e.~$\spr{Au,u}>0$ for all $u\in\operatorname{dom}A$,
$u\neq0$) selfadjoint elliptic $A\in\Psi^m$ with $m>0$, one defines $A^z$ as
$(A^{-1})^{-z}$. It~is then a classical result of Seeley that in both cases,
$A^z$~is a \psdo/ of order~$mz$, with principal symbol~$\sigma(A)^z$.
In~particular, if~we define the space $\HH_A$ as the completion of
$C^\infty(\pOm)$ with respect to the~norm
$$ \|u\|^2_A := \sprb{Au,u} = \|A^{1/2}u\|_\pOm^2,   $$ 
then $\HH_A=W^{m/2}(\pOm)$ as~sets, with equivalent norms.
All~this remains in force also for operators of order $m=0$; note that the
positivity of $A$ then implies, in~particular, that $A$ is injective and,
hence, with bounded inverse on~$L^2(\pOm)$, so~one can again define $A^z$ for
any $z\in\CC$ by the spectral theorem for (bounded) selfadjoint operators.

For $P\in\Psi^m$, the \gto/ $T_P:W^m\hol(\pOm)\to\hdva$ is
defined~as 
$$ T_P = \Pi P,  $$
where $\Pi:L^2(\pOm)\to\hdva$ is the orthogonal projection
(the~Szeg\"o~projection). Alternatively, one~may view $T_P$ as the operator
$$ T_P = \Pi P \Pi   $$
on~all of~$W^m(\pOm)$. Then~$T_P$ maps continuously $W^s(\pOm)$ into 
$W^{s-m}\hol(\pOm)$, for each $s\in\RR$. The~microlocal structure of \gto/s
was described by Boutet de Monvel and Guillemin~\cite{\BdMi}~\cite{\BdMG},
who~proved in particular the following facts. Let~$\Sigma$ denote the
half-line bundle 
$$ \Sigma := \{(x,\xi)\in T^*(\pOm): \xi=t\eta_x, t>0 \} ,  \tag\tSI  $$
where $\eta$ is the restriction to $\pOm$ of the 1-form $\operatorname{Im}
(-\partial\rho)=(\dbar\rho-\partial\rho)/(2i)$; the~strict pseudoconvexity
of~$\Omega$ implies that $\Sigma$ is a symplectic submanifold of the cotangent
bundle~$T^*(\pOm)$. 

\roster
\item"(P1)" For any $T_P$, $P\in\Psi^m$, there in fact exists
$Q\in\Psi^m$ such that $T_P=T_Q$ and $Q$ commutes with~$\Pi$.
\item" "
(Hence, $T_P=T_Q$ is just the restriction of $Q$ to the Hardy space. 
It~follows, in~particular, that generalized Toeplitz operators~$T_P$
form an algebra.)

\item"(P2)" It~can happen that $T_P=T_Q$ for two different \psdo/s $P$ and~$Q$.
However, one~can define unambiguously the order of $T_Q$ as $\min\{\ord(P):\;
T_P=T_Q\}$, and the symbol of $T_Q$ as $\sigma(T_Q):=\sigma(Q)|_\Sigma$
if $\ord(Q)=\ord(T_Q)$.

\item"(P3)" The order and the symbol obey the usual laws: $\ord(T_QT_P)=
\ord(T_Q)+\ord(T_P)$ and $\sigma(T_Q T_P)=\sigma(T_Q)\sigma(T_P)$.

\item"(P4)" If $\ord(T_P)=m$, then $T_P$ maps $W^s\hol(\pOm)$ continuously into
$W^{s-m}\hol(\pOm)$, for any $s\in\RR$. In~particular, if $\ord(T_P)=0$ then
$T_P$ is a bounded operator on~$L^2(\pOm)$; if~$\ord(T_P)<0$, then it is
even compact. 

\item"(P5)" If $P\in\Psi^m$ and $\sigma(T_P)=0$, then there exists
$Q\in\Psi^{m-1}$ with $T_Q=T_P$.

\item"(P6)" We~will say that a generalized Toeplitz operator~$T_P$
is~elliptic if $\sigma(T_P)$ does not vanish.
Then $T_P$ has a parametrix, i.e.~there exists an elliptic generalized Toeplitz
operator $T_Q$ of order $-\ord(T_P)$, with $\sigma(T_Q)=\sigma(T_P)^{-1}$,
such that $T_P T_Q-I$ and $T_Q T_P-I$ are smoothing operators
(i.e.~have Schwartz kernel in $C^\infty(\pOm\times\pOm)$).
\endroster

Note that from (P3) and (P5) we obtain, in~particular, that
$$ \ord[T_P,T_Q] \le \ord(T_P)+\ord(T_Q)-1.    \tag\tCM  $$

For~an elliptic \gto/ $T_P$ of order $m>0$ or~$m<0$ which is positive
selfadjoint as an operator on~$\hdva$, it~again follows from (P1),
(P6) and the result of Seeley recalled above that the complex powers~$T_P^z$,
$z\in\CC$, defined by the spectral theorem, are~elliptic \gto/s of order~$mz$,
with symbol~$\sigma(T_P)^z$ (see~Proposition~16 in~\cite{\Ejfa}
for the details); and, likewise, the~space $\HH_{T_P}$ defined as 
the completion of $C^\infty\hol(\pOm)$ with respect to the~norm
$$ \|u\|^2_{T_P} := \sprb{T_P u,u} = \|T_P^{1/2}u\|_\pOm^2  \tag\tVB   $$
coincides with $W^{m/2}\hol(\pOm)$, with equivalent norms.
The~corresponding space 
$$ \KK\HH_{T_P}:=\{\KK u:u\in\HH_{T_P}\}  $$
of~holomorphic functions on~$\Omega$ thus coincides~with
$$ \KK\HH_{T_P} = \KK W^{m/2}\hol(\pOm) = W^{(m+1)/2}\hol(\Omega),  \tag\tVF $$
with equivalent norms.

We~conclude this section with a simple criterion for Schatten class membership
of~\gto/s.

\proclaim{Proposition~\prodef\PC} A~\gto/ $T_Q$ of order $-q$ on~$\pOm$, $q>0$,
belongs to $\SS^p$ for all $p>n/q$, $n=\dim_\CC\Omega$. \endproclaim

\demo{Proof} Choose a positive selfadjoint elliptic \gto/ of order $-1$
on~$\pOm$ with positive symbol, for~instance, $T_\Lambda$~where
$\Lambda=\KK^*\KK$, cf.~the beginning of the next section.
Then $T_\Lambda^{-q}T_Q$ is a bounded operator;
since $\SS^p$ is an ideal, it~therefore suffices to show that
$T_\Lambda^q\in\SS^p$ for $p$ as indicated.

To~prove the latter, we~proceed as in Theorem~3 in \cite{\EZ}: namely,
let $0<\lambda_1\le\lambda_2\le\dots$ be the eigenvalues of~$T_\Lambda^{-1}$,
counting multiplicities, and denote
$$ N(\lambda) = \operatorname{card} \{j:\,\lambda_j<\lambda \}  $$
the~corresponding counting function. By~Theorem~13.1 in~\cite{\BdMG},
$$ N(\lambda) = c \lambda^n + O(\lambda^{n-1})
 \qquad \text{as } \lambda\to+\infty    $$
with some positive constant~$c$; implying that
$$ \frac1\lambda = \Big(\frac{N(\lambda)}c\Big)^{-1/n}
 [1+O(N(\lambda)^{-1/n})] .  $$
Consequently,
$$ \align
\|T_\Lambda^q\|_{\SS^p}^p &= \sum_{j=1}^\infty \lambda_j^{-pq}   
 = \int_{[\lambda_1,\infty)} \lambda^{-pq} \,dN(\lambda)   \\    
&= \int_1^\infty \Big(\frac cN\Big)^{pq/n} [1+O(N^{-1/n})] \, dN, \endalign  $$
which is finite for $\frac{pq}n>1$, i.e.~for $p>\frac nq$.   \qed  \enddemo

From the last proof one can in fact show that $T_\Lambda^q$ and, hence,
$T_Q$~belongs to the ideal $\SS^{n/q,\infty}$ of operators $T$ whose singular
numbers satisfy $s_j(T)=O(j^{-q/n})$ as $j\to\infty$, and which is properly
contained in $\SS^p$ for all $p>\frac nq$.

\head 3. Sobolev-Bergman spaces and restrictions\endhead
The~Poisson operator $\KK$ is in particular bounded from $L^2(\pOm)$
into~$L^2(\Omega)$, and we denote by~$\KK^*:L^2(\Omega)\to L^2(\pOm)$
its adjoint. Operators of the~form
$$ \Lambda_w := \KK^* w \KK,   $$
where $w$ is a function on~$\Omega$, are~governed by a calculus developed
by Boutet de Monvel~\cite{\BdMA}. Namely, for~$w$ of the~form
$$ w = \rho^\alpha g, \quad \alpha>-1, \; g\in C^\infty(\oOm),   $$
$\Lambda_w$ is an operator in~$\Psi^{-\alpha-1}$, with principal symbol
$$ \sigma(\Lambda_w)(x,\xi) = \frac{\Gamma(\alpha+1)}{2|\xi|^{\alpha+1}}
 g(x) \|\eta_x\|^\alpha .   $$ 
In~particular, we~obtain that $\Lambda:=\Lambda_{\jedna}=\KK^*\KK$ is~an
elliptic operator in~$\Psi^{-1}$, and more generally, $\Lambda_{\rho^\alpha}
=\KK^*\rho^\alpha\KK$ is~an elliptic operator in~$\Psi^{-\alpha-1}$,
for~any $\alpha>-1$.
From the simple computation
$$ \aligned
\int_\Omega |\KK u|^2 \, w \, dz &= \spro{w\KK u,\KK u}  \\
&= \sprb{\Lambda_w u,u}  \\
&= \sprb{T_{\Lambda_w}u,u},   \endaligned   \tag\tEE   $$
valid for any $u\in C^\infty\hol(\pOm)$, and~(\tVF) we thus see that the space
$$ A^2_{\alpha,\rho}(\Omega) := L^2\hol(\Omega,\rho^\alpha\,dz),
 \qquad \alpha>-1,   $$
coincides with $W^{-\alpha/2}\hol(\Omega)$, with equivalent norms,
independently of the choice of the defining function~$\rho$. 
This~suggests extending the definition of the spaces $A^2_{\alpha\rho}$ in
this manner to all real~$\alpha$: namely, let~us introduce the notation
$$ \adah := W^{-\alpha/2}\hol(\Omega),  \qquad \alpha\in\RR.  $$ 
It~was shown by Beatrous~\cite{\Bea} for smoothly bounded strictly pseudoconvex
domains $\Omega$ in a Stein manifold $W$ that there exist many equivalent norms
on~$\adah$, beside the Sobolev norm inherited from~$W^{-\alpha/2}(\Omega)$.
Namely, if~$m$ is a nonnegative integer and
$m>-\frac{\alpha+1}2$, then $f\in\adah$ if~and only~if $\partial^\nu f$
belongs to $A^2_{\alpha+2m,\rho}$ for all multiindices $\nu$ with $|\nu|\le m$,
and the norm in $\adah$ is equivalent~to
$$ \|f\|_{\alpha\#m\rho} := \Big(\sum_{|\nu|\le m}
 \|\partial^\nu f\|_{\alpha+2m,\rho} ^2 \Big) ^{1/2}.   \tag\tWC   $$
Furthermore, in~fact one need not consider all the derivatives in~(\tWC),
but~only ``radial''~ones: namely, if~$\cD$ is the holomorphic vector field
on~$\Omega$ given~by
$$ \cD := \sum_{j=1}^n (\dbar_j\rho) \partial_j ,   $$ 
then a holomorphic function $f$ belongs to $\adah$ if and only if $\cD^j f\in
L^2_{\alpha+2m,\rho}$ for all $0\le j\le m$, and
$$ \|f\|_{\alpha\flat m\rho} := \Big( \sum_{j=0}^m \|\cD^jf\|
 _{\alpha+2m,\rho}^2 \Big)^{1/2}  \tag\tWE   $$
is~an equivalent norm in~$\adah$. Here $L^2_{\alpha\rho}:=L^2(\Omega,
\rho^\alpha dz)$ for $\alpha>-1$, and~we use $\|\cdot\|_{\alpha\rho}$
to~denote the norm in~$L^2_{\alpha\rho}$. Again, both~in (\tWC) and in~(\tWE),
$\rho$~can be an arbitrary defining function, and different choices of
$\rho$ lead to equivalent norms. We~remark that a proof of all the above
facts can be given based on (\tVF) and the machinery of \gto/s reviewed in
the preceding section (which~is completely different from the methods used
in~\cite{\Bea}): namely, one~checks that the norms in (\tWC) and~(\tWE) are
special cases of the norm~(\tVB), with
$$ P = P_{\alpha\#m\rho} := \sum_{|\nu|\le m} \KK^*
 \dbar{}^\nu \rho^{\alpha+2m} \partial^\nu \KK     \tag\tWI   $$
and
$$ P = P_{\alpha\flat m\rho} := \sum_{j=0}^m \KK^*
 \cD^{*j} \rho^{\alpha+2m} \cD^j \KK ,   \tag\tWJ   $$
respectively, and that $P$ is a positive selfadjoint elliptic \psdo/ on $\pOm$
of order $-\alpha-1$; see~Sections~5--7 in~\cite{\Ejfa} for the details.

Finally, \cite{\Bea} also gives a result concerning restrictions of functions
in $\adah$ to complex submanifolds which intersect $\pOm$ transversally. 
(In~fact \cite{\Bea} treats even the case of $L^p$-Sobolev spaces of
holomorphic functions for any $p>0$, not only $p=2$.)
Namely, if~$V$ is such a submanifold in a neighbourhood of~$\oOm$,
then Corollary~1.7 in \cite{\Bea} asserts
that the restriction map
$$ R_V : f\longmapsto f|_V   \tag\tTV   $$
actually sends each $\adah(\Omega)$ continuously onto~$\akov$: 
$$ R_V \adah(\Omega) = \akov \quad\text{continuously,}
 \qquad \forall \alpha\in\RR,   \tag\tTR   $$
where
$$ k = n - \dim_\CC V   $$ 
is~the codimension of~$V$ in the $n$-dimensional Stein manifold~$W$.

We~will need a somewhat more precise information on the nature of the
restriction operator
$R_V$ and its relationships to the inner products like (\tAC), (\tWC), (\tWE)
on $\adah(\Omega)$ and~$\akov$. To~that~end, we~now~review some properties of
the Szeg\"o projection $\Pi:L^2(\pOm)\to\hdva$ due to Boutet de Monvel and
Sj\"ostrand~\cite{\BdMS}. 

Recall that a \emph{Fourier integral distribution} is~an integral of the form
$$ u(x) = I(a,\phi)(x) := \int e^{i\phi(x,\theta)} a(x,\theta) \,d\theta.
 \tag\tFA   $$
Here $a$ is a (classical) symbol in the H\"ormander class $S^m(U\times\RR^N)$,
$U\subset\RR^n$, $m\in\RR$, and $\phi\in C^\infty(U\times\dR^N)$, $\dR^N:=\RR^N
\setminus\{0\}$, is~a~\emph{nondegenerate phase function}, meaning that $\phi$
is real-valued, $\phi(x,\lambda\theta)=\lambda\phi(x,\theta)$ for $\lambda>0$,
$d_{(x,\theta)}\phi\neq0$, and $d_{(x,\theta)}\frac{\partial\phi}{\partial
\theta_j}$, $j=1,\dots,N$, are~linearly independent on the set where
$d_\theta\phi=0$. The~integral (\tFA) converges absolutely when $m<-N$,
and can be defined as a distribution on $U$ for any real~$m$.
The~image $\Lambda_\phi$ of the set $\{(x,\theta):d_\theta\phi(x,\theta)=0\}$
under the map $(x,\theta)\mapsto(x,d_x\phi(x,\theta))$ is then a conical
Lagrangian submanifold of $\Tb U:=T^*U\setminus\{0\}$, the~cotangent bundle
of $U$ with zero section removed. The~set of all distributions of the form
(\tFA) turns out to depend not on $\phi$ but only on~$\Lambda_\phi$, modulo
smooth functions: namely, if $\psi\in C^\infty(U\times\dR^M)$ is another phase
function such that $\Lambda_\phi=\Lambda_\psi$ in a neighbourhood of
$(x_0,\xi_0)\in\Tb U$, and $a\in S^m(U\times\RR^N)$ is supported in a small
conical neighbourhood of~$(x_0,\xi_0)$, then there exists
$b\in S^{m'}(U\times\RR^M)$, where $m'=m+(M-N)/2$, supported in a small conical
neighbourhood of $(x_0,\xi_0)$, such that $I(a,\phi)-I(b,\psi)\in C^\infty(U)$.
Moreover, if $a$ is elliptic, then so~is~$b$. Given a conical Lagrangian
submanifold $\Lambda$ of~$\Tb U$, one~can therefore unambiguously define
the space of associated Fourier integral distributions
$$ I^m(U,\Lambda) := \{u=I(a,\phi)+v\text{ locally, }\Lambda_\phi=\Lambda, \,
 a\in S^{m-\frac N2+\frac n4}, \, v\in C^\infty \},   \tag\tFB   $$
and its subset $I^m\ell$ with $a$ elliptic.
(The~reason for the shift by $\frac n4$ will become apparent in a moment.)
The~whole construction carries over in a straightforward manner from subsets
$U\subset\RR^n$ also to real manifolds of dimension~$n$.

If~$X,Y$ are two compact real manifolds, and $\Lambda$ is a conical Lagrangian
submanifold of $\Tb X\times\Tb Y\subset\Tb(X\times Y)$, an~operator from
$C^\infty(X)$ into $\cD'(Y)$ whose distributional (Schwartz) kernel belongs
to $I^m(Y\times X,\Lambda)$ is~called a Fourier integral operator
(FIO~for short) of~order~$m$. The~set
$$ C := \{(x,\xi,y,-\zeta): \; ((x,y),(\xi,\zeta))\in\Lambda \}
 \subset \Tb X \times \Tb Y   $$
is~called the ``canonical relation'' corresponding to~$\Lambda$, and~we denote
the space of all FIOs (with classical symbols) from $X$ into $Y$ of order $m$
and with canonical relation $C$ by $I^m(X,Y,C)$, and by $I^m\ell(X,Y,C)$ its
subset of elliptic elements. The~$L^2$ adjoint (with respect to some smooth
volume elements on $X$ and~$Y$) of~an~operator $A\in I^m(X,Y,C)$ belongs
to~$I^m(Y,X,C^t)$, where $C^t=\{(y,\zeta,x,\xi): (x,\xi,y,\zeta)\in C\}$.
If~$X,Y,Z$ are three compact real manifolds and $A_1\in I^{m_1}(X,Y,C_1)$,
$A_2\in I^{m_2}(Y,Z,C_2)$, where $C_1$ and~$C_2$ intersect nicely%
\footnote{In~detail: if~$\Delta$ denotes the ``diagonal''
$\Delta=\{(a,b,b,c): a\in\Tb X, b\in\Tb Y, c\in\Tb Z\}$, then
(i) $C_1\times C_2$ should intersect $\Delta$ transversally
(i.e.~the~sum of the tangent spaces should be equal to the full tangent space
of $\Tb X\times(\Tb Y)^2\times\Tb Z$ at each point of intersection), and (ii)
the~natural projection $(C_1\times C_2) \cap\Delta\to\Tb(X\times Z)$ should be
injective and proper; its~image is denoted~$C_1\circ C_2$.}, then
$$ A_2 A_1 \in I^{m_1+m_2}(X,Z,C_1\circ C_2),   \tag\tFC   $$
and similarly for $I^m$ replaced by~$I^m\ell$.
(This composition~law ---~i.e.~that $I^{m_1}\circ I^{m_2}\subset I^{m_1+m_2}$
--- is~the reason for the shift by $\frac n4$ in~(\tFB).)

Pseudodifferential operators are special case of FIOs corresponding to the
phase function $\phi(x,y,\theta)=(x-y)\cdot\theta$; thus $C=\diag\Tb X$ is
the identity relation $\{(x,\xi,x,\xi): (x,\xi)\in\Tb X\}$, and \psdo/s can
be composed with any~FIO, yielding a FIO with the same canonical relation as
the original~FIO.

Finally, up~to a number of technicalities which we will not go into here
(see~the references mentioned below for the details), the~calculus of FIOs
extends also to complex valued phase functions~$\phi$ with $\Im\phi\ge0$.
The~technicalities stem from the fact that the set $\{(x,\theta):
d_\theta\phi(x,\theta)=0\}$ is~no longer a (real) manifold in $U\times\RR^N$
in~general, and needs to be replaced, roughly speaking, by~the ``real part''
of~its ``almost analytic'' complex extension; the~same applies to the conical
Lagrangian manifolds $\Lambda_\phi$ and canonical relations~$C$. With~these
modifications, the~whole formalism of Fourier integral distributions and FIOs
just described remains in force also for complex-valued phase functions.

The~reader is referred e.g.~to Grigis and Sj\"ostrand~\cite{\GrSj},
H\"ormander~\cite{\Horm} (Chapter~25), Melin~and Sj\"ostrand~\cite{\MeSj}
and Treves~\cite{\Trev} (Chapters~VIII and~X) for full accounts of the theory
of FIOs with real as well as complex valued phase functions.

The~main result of \cite{\BdMS} then says that for any smoothly-bounded
strictly pseudo\-convex domain $\Omega$ as in Section~2, the~Szeg\"o kernel
$S(x,y)$ is a Fourier integral distribution in
$I^0(\pOm\times\pOm,\diag\Sigma_\pOm)$,
and the Szeg\"o projection $\Pi:L^2(\pOm)\to\hdva$ is an elliptic FIO with
complex valued phase function in $I^0\ell(\pOm,\pOm,\diag\Sigma_\pOm)$, where
$$ \aligned
\diag\Sigma_\pOm &= \{(\Upsilon,\Upsilon): \Upsilon\in\Sigma_\pOm\}  \\
&= \{(x,t\eta_x,x,t\eta_x):\,x\in\pOm,t>0\}\subset \Tb(\pOm)\times\Tb(\pOm),
\endaligned   $$ 
with $\Sigma=\Sigma_\pOm$ as in~(\tSI). More specifically, one~has%
\footnote{For $\Omega=\BB^d$ with $\rho(x,y)=1-\spr{x,y}$, one has simply
$a(x,y,\theta)=\theta^{d-1}/\lambda(\partial\BB^d)$.}
$$ S(x,y) = \int_0^\infty e^{-\theta\rho(x,y)} a(x,y,\theta) \,d\theta,
 \tag\tFE   $$
where $a$ is an elliptic symbol in $S^{n-1}(\oOm\times\oOm\times\RR_+)$
and $\rho(x,y)$ is an ``almost analytic'' extension of the defining
function~$\rho$, namely $\rho(\cdot,\cdot)\in C^\infty(\oOm\times\oOm)$
satisfies $\rho(x,x)=\rho(x)$, $\rho(y,x)=\overline{\rho(x,y)}$,
while $\dbar_x\rho(x,y),\partial_y\rho(x,y)$ vanish to infinite order on the
diagonal $x=y$, and $2\Re\rho(x,y)\ge\rho(x)+\rho(y)+c|x-y|^2$
for all $x,y\in\oOm$ for some $c>0$. It~follows that an (elliptic) \gto/ $T_P$
on $\pOm$ of order $m$ is an (elliptic) FIO in
$I^m(\pOm,\pOm,\diag\Sigma_\pOm)$, and in fact \gto/s of order $m$ on $\pOm$
are precisely those operators $A\in I^m(\pOm,\pOm,\diag\Sigma_\pOm)$ for which
$A=\Pi A\Pi$ ($=\Pi A=A\Pi$); see~\cite{\BdMT},~p.~21,~\S v.
\footnote{Page~253 in~\cite{\BdMi} gives a construction of an operator $H$
from $L^2(\RR^n)$ onto the Hardy space such that $H^*H=I$ while $HH^*=\Pi$;
$H$~is~in fact a FIO of order~0. An~operator $T$ satisfying $T=\Pi T\Pi$
then equals $T=HQH^*$ where $Q=H^*TH$. Now~the paragraph before (1.7) on the
same page 253 of~\cite{\BdMi} outlines a proof that any such $Q$ can be
obtained as $Q=H^*PH$ for some \psdo/~$P$. It~follows that $T=\Pi P\Pi=T_P$
is~a~\gto/, as~claimed.}

After these preparations, we~can state an observation which in some sense is
our main result of this section. Let~$V$ be a complex submanifold in a
neighbourhood of~$\oOm$ which intersects $\pOm$ transversally, $k=n-\dim_\CC V$
and let $R_V$ be the restriction operator~(\tTV). We~will denote the Szeg\"o
projections on $\pOm$ and~$\partial(\Omega\cap V)=\pOcv$ by $\Pi$ and~$\Pi_V$,
respectively, and similarly by $\rho$ and $\rho_V$ the respective defining
functions as well as their almost-analytic extensions (thanks to the
transversality hypothesis, one~can take $\rho_V=\rho|_{V\times V}$,
and~we will assume this from now~on),
by~$\eta=\frac1{2i}(\dbar\rho-\partial\rho) |_{\pOm}$ and 
$\eta^V=\frac1{2i}(\dbar\rho_V-\partial\rho_V) |_{\pOcv}$
the corresponding one-forms on the boundary,
$\KK$~and $\KK_V$ the respective Poisson operators, etc.
Finally, we~denote~by
$$ \rdv := \gamma_V R_V \KK : u \longmapsto u|_\pOcv   $$
the~action of $R_V$ on boundary values.

\proclaim{Proposition~\prodef\PR} $\rdvpb=\Pi_V\rdvpb$~is an elliptic FIO
from $\pOm$ to $\pOcv$ of order $k/2$, with canonical relation\footnote
{More precisely, (\prtFF)~is the ``real part'' of the canonical relation.} 
$$ \sibo := \{(x,t\eta_x,x,t\eta^V_x): x\in\pOcv, t>0 \}
 \subset \Tb(\pOm)\times\Tb(\pOcv)   \tag\tFF   $$
$($the~``restriction'' of $\diag\Sigma_{\pOm}$ to
$\Tb(\pOm)\times\Tb(\pOcv))$.        \chk\prtFF\tFF

Furthermore, if~$T$ is a \gto/ on $\pOm$ of order $s\in\RR$, then
$(\rdvpb)T(\rdvpb)^*$ is a \gto/ on $\pOcv$ of order~$s+k$,
which is elliptic if~$T$~is.%
\footnote{Throughout the paper, unless explicitly stated otherwise,
the~adjoint $X^*$ of an operator $X$ acting on a Hilbert space of functions
on a domain $\Omega$ or its boundary~$\pOm$ (or~between two such spaces)
is~always meant with respect to the $L^2$ products on $\Omega$ or~$\pOm$.
This is in line with the standard convention in \psdo/ theory --- if~$X$ is
a~\psdo/ of order~$m$, then $X^*$ is also of order~$m$, and both $X$ and $X^*$
(sic!) map $W^s$ into~$W^{s-m}$ for any real~$s$. (Normally, the~adjoint of
$X:W^s\to W^{s-m}$ would be $X^*:W^{s-m}\to W^s$, the~reason of course being
that the latter adjoint is taken with respect to the $W^s$ and $W^{s-m}$ inner
products and not with respect to the $L^2$ products. The~only place where we
will use the genuine, instead of $L^2$, adjoints are the operators $T^*$ in
the proofs of Theorem~{\prPY} in Section~4 and of Main Theorem in Section~5.)
\endgraf
Strictly speaking, by the $L^2$ inner product above we also mean its extension
to the duality pairing between $W^s$ and~$W^{-s}$, $s\in\RR$, which coincides
with the $L^2$ pairing when both arguments are smooth functions.}  \endproclaim

Note that the second part of the proposition would actually follow immediately
by~(\tFC) if the canonical relations $\sibo$, $\sibo^t$ and $\diag\Sigma_\pOm$
intersected nicely. However this does not seem to be the case (unless~$k=0$);
fortunately it is possible to give a direct~proof.

\demo{Proof} Set~temporarily, for brevity, $A:=\rdvpb$. Since the restriction
of a holomorphic function to a complex submanifold is again holomorphic, it~is
clear that $A=\Pi_V A$. Also~by~(\tFE), the~Schwartz kernel of $A$ is simply
the restriction
$$ R_{\partial V,x} S(x,y) = \int_0^\infty e^{-\theta\rho(x,y)}
 a(x,y,\theta) \, d\theta, \qquad x\in\pOcv, y\in\pOm.   $$
Comparing this with (\tFA) shows that, just as~(\tFE), the~right-hand side
is a Fourier integral distribution, and, hence, $A$~is a~FIO, with phase
function $i\theta\rho(x,y)$, $x\in\pOcv$, $y\in\pOm$, and canonical relation
given~by~(\tFF). Since $a$ is an elliptic symbol in~$S^{n-1}$, the~order of
$\rdvpb$~is $n-1+\frac{\dim\RR_+}2-\frac{\dim\pOm+\dim\pOcv}4 =
n-1+\frac12-\frac{(2n-1)+(2n-2k-1)}4 = \frac k2$, proving the first
part of the proposition.

For~the second part, note that from the formulas
$$ Af(x) = \int_{y\in\pOm} f(y) S(x,y) \,d\lambda(y), \qquad
 A^*g(y) = \int_{x\in\pOcv} g(x)S(y,x) \,d\lambda_V(x),   $$
we~get
$$ ATA^* f(x) = \int_{y\in\pOm} S(x,y) T_y \int_{x_1\in\pOcv} f(x_1)
 S(y,x_1) \,d\lambda_V(x_1) \,d\lambda(y),   $$
where the subscript $y$ in $T_y$ refers to the variable $T$ is being
applied~to. Thus~the Schwartz kernel of $ATA^*$~is
$$ \Cal K_{ATA^*}(x,x_1) := \int_{y\in\pOm} S(x,y) T_y S(y,x_1) \,d\lambda(y)
= (TS_{x_1})(x)  $$
by~the reproducing property of the Szeg\"o kernel, where $S_y(x):=S(x,y)$.
Now~we may assume that $T=T_P$ for some \psdo/ $P$ of the same order which
commutes with~$\Pi$; and by the standard symbol calculus for \psdo/s
(see,~for~instance, Theorem~4.2 in H\"ormander~\cite{\Hcpam})
we~have quite generally
$$ T_x \int_0^\infty e^{-\theta\rho(x,y)} a(x,y,\theta) \,d\theta 
 = \int_0^\infty e^{-\theta\rho(x,y)} b(x,y,\theta) \,d\theta   \tag\tSZ  $$
where $b\in S^{m+s}(\oOm\times\oOm\times\RR_+)$ if $a\in
S^m(\oOm\times\oOm\times\RR_+)$, and with $b$ elliptic if $a$ is elliptic.
(See~the proof of Theorem~5 in \cite{\Epay} for the details.)
Restricting to $x,y\in\oOm\cap V$ we thus see that $ATA^*$ is a~FIO,
elliptic if $T$ is elliptic, of order $s+k$ and with the same canonical
relation as~$\Pi_V$, i.e.~$\diag\Sigma_{\pOcv}$. Since $A=\Pi_V A$ and,
hence, $ATA^*=\Pi_V ATA^*\Pi_V$, it~therefore follows (see~the end of
the paragraph after (\tFE)~above) that $ATA^*$ is a \gto/ on~$\pOcv$,
proving the second part of the proposition.   \qed   \enddemo

\proclaim{Corollary~\prodef\PZ} Under the same hypotheses as for the
preceding proposition, $A=\rdvpb$ is bounded from $W^{-s/2}(\pOm)$
into $W^{-(s+k)/2}\hol(\pOcv)$, for any $s\in\RR$, and its range has
finite codimension.

In~particular, $\rdv:W^{-s/2}\hol(\pOm)\to W^{-(s+k)/2}\hol(\pOcv)$
and $R_V: \adah(\Omega)\to\akov$ are bounded for any $s,\alpha\in\RR$,
with ranges of finite codimension.    \endproclaim

\demo{Proof} Choose an invertible elliptic \gto/ $T_P$ of order $\frac s2$
on~$\pOm$ (for~instance, $\Lambda^{-s/2}$ with $\Lambda=\KK^*\KK$ as before),
and similarly an invertible elliptic \gto/ $T_Q$ of order $-\frac{k+s}2$
on~$\pOcv$. Then by the last proposition and the properties of \gto/s,
$T_QAT_PT^*_PA^*T^*_Q$ is a~\gto/ of order~0, hence, a~bounded operator
on~$H^2(\pOcv)$. Since an operator $X$ between Hilbert spaces is bounded
if and only if $XX^*$~is, $T_QAT_P$~must be bounded from $H^2(\pOm)$
into~$H^2(\pOcv)$. By~the mapping properties of \gto/s again, this means
that $A$ is bounded from $W^{-s/2}\hol(\pOm)$ into $W^{-(s+k)/2}\hol(\pOcv)$,
proving the first claim. 
Furthermore, since $T_QAT_PT^*_PA^*T^*_Q$ is elliptic, by~the property
(P6) of \gto/s it has a parametrix and, hence, is~a~Fredholm operator
on~$H^2(\pOcv)$; thus it has (closed) range of finite codimension. 
Since again an operator $X$ has range of finite codimension if and only
if $XX^*$~does, while $T_P$ and $T_Q$ are isomorphisms of $H^2(\pOm)$ onto
$W^{-s/2}\hol(\pOm)$ and of $W^{-(s+k)/2}\hol(\pOcv)$ onto $H^2(\pOcv)$,
respectively, the~second claim about $A$ follows. The~second half of the
corollary is immediate from the first.   \qed   \enddemo

From the results of Beatrous, we~know that for $\Omega$ in $\CC^n$ or in
a~Stein manifold, the~restriction operator $R_V:\adah(\Omega)\to\akov$ is
actually~onto. On~the other hand, consider the situation when $W$ is the
tautological line bundle over the complex projective space~$\CC P^1$,
i.e.~$W=\{(\CC z,cz): c\in\CC, z\in\CC^2, z\neq0\}$, let $\Omega$ be the
unit disc bundle $\{(\CC z,cz)\in W:\|cz\|<1\}$, and take $V=\{(\CC z,cz)
\in W: c\in\CC, z_1^2-z_2^2=0\}$. Then $\Omega\cap V$ consists of the two
fibers of $\Omega$ over the points $(1:1)$ and $(1:-1)$ of~$\CC P^1$,
i.e.~two~disjoint discs. The~function equal to 0 on one disc and to 1 on
the other one is holomorphic in~$\Omega\cap V$, but~cannot be the restriction
to $V$ of a holomorphic function $f$ on~$\Omega$: any~such $f$ must be constant
on the zero section of~$\Omega$ (which is a compact complex submanifold
of~$\Omega$), hence assumes one and the same value in the centers of the
two discs that form $\Omega\cap V$. Thus finite codimension of the range of
$R_V$ is indeed the best one can get (in~this example, the~codimension is~1).

The~statement of the last corollary should be contrasted with the situation
for full Sobolev spaces~$W^s$ (instead of~$W^s\hol$): there the restriction
map $R_V$ maps $W^s(\Omega)$ into $W^{s-(k/2)}(\Omega\cap V)$ (and~$\rdv$
maps $W^s(\pOm)$ into $W^{s-(k/2)}(\pOcv)$) only for $s>k/2$, by~the Sobolev
trace theorem; whereas for the subspaces $W^s\hol$ of holomorphic functions
this holds for all real~$s$. This is completely parallel to what happens
for the ``boundary-value'' operator~$\gamma$ from Section~2, which maps
$W^s(\Omega)\to W^{s-1/2}(\pOm)$ only for $s>1/2$, but $W^s\hol(\Omega)
\to W^{s-1/2}\hol(\pOm)$ and $W^s\harm(\Omega)\to W^{s-1/2}(\pOm)$
for any real~$s$.

Finally, note that when $R_V$ is onto, then $R_V R^*_V$ must be invertible
(by~Banach's inverse mapping theorem), and $R^*_V(R_V R^*_V)^{-1}$, being
a~right inverse to~$R_V$, is~then a bounded extension operator from $\akov$
into~$\adah(\Omega)$, for any $\alpha\in\RR$.

\head 4. The~case of a smooth submanifold\endhead
Throughout the rest of this paper, the~space $\adac$ on~$\BB^d$, 
$\alpha\in\RR$, will always be understood to be equipped with a~norm 
of the form
$$ \|f\|^2 = \sprdb{T_Y\gamma f,\gamma f}    \tag\tQQ   $$
for some positive selfadjoint elliptic \gto/ $T_Y$ on~$\partial\BB^d$
of order~$-\alpha-1$. In~particular, this includes the weighted Bergman
spaces $L^2\hol(\BB^d,\rho^\alpha\,dz)$ for $\alpha>-1$,
with~$Y=\Lambda_{\rho^\alpha}$ (cf.~(\tEE)),
the~Hardy space on~$\partial\BB^d$ (with~$Y$ the identity operator),
or~any of the Sobolev norms (\tWC) or~(\tWE), for any $\alpha\in\RR$
(with $Y$ given by~(\tWI) and~(\tWJ), respectively). It~also includes
the original norms $\|f\|_{\alpha\circ}$ from~(\tAC): namely, by~computing
the inner products $\spr{z^\nu,z^\mu}$, where $\mu,\nu$ are two multiindices,
in~the Hardy space $H^2(\partial\BB^d)$ and in the unweighted Bergman
space $L^2\hol(\BB^d)$ and comparing the results, one~sees that~$T_\Lambda$,
$\Lambda=\KK^*\KK$, is~the operator on $H^2(\partial\BB^d)$ diagonalized
by the monomial basis $\{z^\nu\}_\nu$ with eigenvalues
$$ T_\Lambda z^\nu = \frac1{2(|\nu|+d)} \; z^\nu .   $$
Thus if $F\in C^\infty(\RR)$ is a function satisfying
$$ F(2x+2d)=\frac{\Gamma(x+d)}{\Gamma(x+1)(x+1)^{d+\alpha}}
 \quad\text{for }x\ge0, \qquad F(x)=0\quad\text{for }x\le-1,   $$
then the operator $B:=F(T_\Lambda^{-1})$ (defined by the functional calculus
for selfadjoint operators) satisfies
$$ \sprdb{Bz^\nu,z^\mu} = \spr{z^\nu,z^\mu}_{\alpha\circ} .   $$
Thus the norm (\tAC) is of the form (\tQQ) with $B$ in the place of~$T_Y$.
Now~by the property (P1) of \gto/s, there exists an elliptic \psdo/ $Q$ of
order $-1$, commuting with~$\Pi$, such that $T_Q=T_\Lambda$. By~elementary
properties of the functional calculus and since $\Pi Q=Q\Pi$, this implies
$B=F(T_Q^{-1})=T_{F(Q^{-1})}$. Finally, the~function $F$ is easily seen to 
belong to the H\"ormander class $S^{-\alpha-1}(\RR)$, so~by the classical
result of Strichartz~\cite{\Strich}, $F(Q^{-1})$~ is also an elliptic \psdo/,
of~order $-\alpha-1$. So~taking $Y=F(Q^{-1})$ shows that the norm (\tAC) is
of the form~(\tQQ), as~claimed. (Note~that for $\alpha=-d$ this includes,
in~particular, also the Drury-Arveson norm $\|f\|_{DA}$ that interests us
most of~all.) A~slight modification of this construction (with $F(2x+2d)=
\Gamma(x+d)/\Gamma(x+d+\alpha+1)$ for $x\ge0$) likewise shows that (\tQQ)
includes also the norms in $\ada(\BB^d)$ for all $\alpha>-d-1$.

The~proof of the result below uses the same ideas as that of our main theorem
in the next section, but~is somewhat simpler.

\proclaim{Theorem~\prodef\PY} Let $V$ be a complex submanifold of $\CC^d$ of
dimension~$n$ that intersects $\partial\BB^d$ transversally, $\alpha\in\RR$,
and $\MM$ the subspace in $\adac$ of functions vanishing on~$V\cap\BB^d$.
Then $[S_j,S^*_k]\in\SS^q$, $j,k=1,\dots,d$, for~all $q>n$.
\endproclaim   \chk\prPY\PY

\demo{Proof} Let $R_V:f\mapsto f|_V$ be the restriction map (\tTV) for
$\Omega=\BB^d$. Then $\MM=\Ker R_V$, and by the result of Beatrous 
in~\cite{\Bea}, we~know that $R_V$ maps $\adac$ boundedly onto~$\akbv$,
$k=d-n$. The~restriction of $R_V$ to $\MMp=(\Ker R_V)^\perp$ is thus an
injective map of $\MMp$ onto~$\akbv$.
Keeping the notations from the end of Section~3, let~$T_X$ be a positive
selfadjoint elliptic \gto/ of order $-\alpha-k-1$ on $\partial(\BB^d\cap V)=
\pbcv$ so that $\spr{T_X\gamma_V f,\gamma_V g}_\pbcv$ is an equivalent inner
product in~$\akbv$ (for~instance, $T_X$~can be one of the operators (\tWI) or
(\tWJ) corresponding to the inner products (\tWC) and~(\tWE), respectively). 
Similarly, as~discussed at the beginning of this section, let~$T_Y$ be a 
positive selfadjoint elliptic \gto/ of order $-\alpha-1$ on $\partial\BB^d$
such that $\sprdb{T_Y\gamma f,\gamma g}$ is the inner product
in~$\adac$. The~composed~map
$$ T := T_X^{1/2} \gamma_V R_V = T_X^{1/2} \rdv \gamma: \adac\to H^2(\pbcv)  $$
then satisfies $\Ker T=\MM$, maps boundedly $\adac$ onto~$H^2(\pbcv)$ and is
an isomorphism of $\MMp$ onto~$H^2(\pbcv)$. Let~$T^*:H^2(\pbcv)\to\adac$ be
its adjoint. By~abstract operator theory, $TT^*$~is then invertible and for
all $f,g\in\adac$,
$$ \spr{P_\MMp f,g}_{\adac} = \spr{(TT^*)^{-1}Tf,Tg}_\pbcv.   \tag\tFG   $$
(Indeed, the~restriction $\tau$ of $T$ to $\MMp=(\Ker T)^\perp$ is injective
and onto, hence invertible; as~$TT^*=\tau\tau^*$, the~invertibility of $TT^*$
follows. As~for~(\tFG), both sides vanish if $f$ or $g$ belongs to~$\MM=
\Ker T$; while for $f,g\in\MMp$, the~right-hand side coincides with
$\spr{(\tau\tau^*)^{-1}\tau f,\tau g}=\spr{\tau^{*-1}\tau^{-1}\tau f,\tau g}
=\spr{f,g}$.) We~claim that
$$ \gathered
TT^* = T_X^{1/2} T_Q T_X^{1/2} \quad \text{for an elliptic generalized} \\
\text{Toeplitz operator $T_Q$ on $\pbcv$ of order } \alpha+k+1.
\endgathered   \tag\tFH   $$
To~see this, note that for any $u\in H^2(\pbcv)$ and $f\in\adac$,
$$ \align
\spr{T^*u,f}_{\adac} &= \spr{u,Tf}_\pbcv
 = \spr{u,T_X^{1/2}\rdvpb\gamma f}_\pbcv \\
&= \sprdb{(\rdvpb)^*T_X^{1/2} u,\gamma f} ,  \endalign   $$
while
$$ \spr{T^*u,f}_{\adac} = \sprdb{T_Y\gamma T^*u,\gamma f} $$
by~the definition of the inner product in~$\adac$. Thus
$$ T^* = \KK T_Y^{-1} (\rdvpb)^* T_X^{1/2}  $$
and
$$ TT^* = T_X^{1/2} \rdvpb T_Y^{-1} (\rdvpb)^* T_X^{1/2}.   $$
Since 
$$ T_Q := (\rdvpb) T_Y^{-1} (\rdvpb)^*   $$
is~an elliptic \gto/ of order $\alpha+k+1$ by Proposition~{\PR}, 
(\tFH)~follows. 

For~the compressions $S_j$ of $M_{z_j}$ to~$\MMp$, $j=1,\dots,d$, we~thus
obtain by (\tFG) and (\tFH) for any $f,g\in\MMp$,
$$ \align
\spr{S_j f,g}_\MMp &= \spr{P_\MMp M_{z_j} f,g}_{\adac}  \\
&= \spr{(T_X^{1/2}T_QT_X^{1/2})^{-1}T_X^{1/2}\gamma_V
 R_V M_{z_j}f,T g}_\pbcv  \\
&= \spr{T_X^{-1/2}T_Q^{-1} \gamma_V R_V M_{z_j}f,T g}_\pbcv  \\
&= \spr{T_X^{-1/2}T_Q^{-1} z_j\gamma_V R_V f,T g}_\pbcv  \\
&= \spr{(T_X^{-1/2}T_Q^{-1} z_jT_X^{-1/2})T f,T g}_\pbcv ,  \endalign   $$
where, abusing notation, we~have used $z_j$ to denote also the operator
of multiplication by (the~boundary value~of) the~restriction $z_j|_\pbcv$
on~$\pbcv$. Thus under the isomorphism $T:\MMp\to H^2(\pbcv)$,
$$ S_j = T^* \cT_j T   $$
where
$$ \cT_j := T_X^{-1/2}T_Q^{-1}T_{z_j}T_X^{-1/2}   $$
is~a~\gto/ on $\pbcv$ of order
$$ \frac{\alpha+k+1}2 - (\alpha+k+1) + 0 + \frac{\alpha+k+1}2 = 0.  $$
Hence, for all $j,l=1,\dots,d$, by~an elementary computation
$$ \align
[S_j,S^*_l] &= [T^* \cT_j T, T^* \cT^*_l T]   \\
&= T^* (\cT_j[TT^*,\cT^*_l] + [\cT_j,\cT^*_l]TT^*
 + \cT^*_l[\cT_j,TT^*]) T.  \endalign   $$
As~$\cT_j$ (and, hence,~$\cT^*_l$) are \gto/s of order~0, and so is~$TT^*$
(by~(\tFH)), by~(\tCM) the last three commutators are \gto/s of order~$-1$
(or~less). Since $\SS^p$ is an ideal and $T$ is bounded, taking $q=1$ and
$\Omega=\BB^d\cap V$ in Proposition~{\PC} thus yields
$$ [S_j,S^*_l] \in \SS^p \quad \forall p>n,  $$
proving the theorem.   \qed   \enddemo

\head 5. Proof of Main Theorem\endhead
Let~now $V$ be as in Main Theorem, i.e.~a~homogeneous variety in $\CC^d$ such
that $V\setminus\{0\}$ is a complex submanifold of $\CC^d\setminus\{0\}$ of
dimension~$n$. Our~idea~is, loosely speaking, to~proceed as in the preceding
proof after removing (``blowing~up'') the singularity of $V$ at the origin.

For $z\in\CC^d\setminus\{0\}$, denote by $\CC z=\{cz:c\in\CC\}$ the
one-dimensional complex subspace through~$z$; the~set of all such subspaces
is the complex projective space~$\CC P^{d-1}$. The~hypotheses on $V$ mean
precisely that $\cV:=\{\CC z:0\neq z\in V\}$ is a complex submanifold
of~$\CC P^{d-1}$ (of~dimension~$n-1$). Consider the tautological line bundle
$\cL$ over $\CC P^{d-1}$, i.e.~the~fiber over a point $\CC z\in\CC P^{d-1}$
is the very same complex line $\CC z$; in~other words, $\cL$~consists of all
points $(\CC z,cz)\in\CC P^{d-1}\times\CC^d$ of the form $(\CC z,cz)$,
$z\in\CC^d\setminus\{0\}$, $c\in\CC$. (The~only role of $c$ is to allow the
second coordinate to be also~0.) Let~$\cL_V:=\{(\CC z,cz):z\in V\setminus\{0\},
c\in\CC\}$ be the part of $\cL$ lying over~$\cV$. 
Finally let $\cB:=\{(\CC z,cz)\in\cL: \|cz\|<1\}$ and $\Omega:=\cB\cap\cL_V$
be~the unit disc bundles of $\cL$ and~$\cL_V$, respectively. Then $\cL$ and
$\cL_V$ are complex manifolds of dimensions $d$ and~$n$, respectively, 
the~subsets $\cB\subset\cL$ and $\Omega\subset\cL_V$ are strictly-pseudoconvex
domains with smooth boundary, and $\Omega$ contains no singular points ---
the~origin has been blown up into the zero section of~$\cL_V$. The~domains
$\cB$ and $\Omega$ will play the same roles as $\BB^d$ and $\bcv$ did in the
proof of Theorem~\PY.

The~map $\pi$ sending $(\CC z,cz)$ into $cz$ sends $\Omega$ back into $\bcv$
(and~$\cB$ into~$\BB^d$), and is bijective except for the zero section
which is taken into the origin. This map translates holomorphic functions
on $\Omega$ into holomorphic functions on~$\bcv$ (any~such function must
be constant on the zero section, as~the latter, being a closed submanifold
of~$\CC P^{d-1}$, is~a~compact complex manifold, so~the translated function
on $\bcv$ is single-valued at the origin). Similarly, $\cB$~is mapped
to~$\BB^d$, bijectively except for the zero section being mapped into
the origin, and holomorphic functions on $\cB$ correspond precisely to
holomorphic functions on~$\BB^d$.

One~defines harmonic functions on~$\cB$ as those annihilated by
$\dbar{}^*\dbar$, where $\dbar{}^*$ is the formal adjoint of $\dbar$
with respect to some Hermitian metric on~$\cL$, for~instance the Cartesian
product of the Fubini-Study metric on $\CC P^{d-1}$ and the Euclidean metric
in the fibers; similarly for~$\Omega$. It~should be noted that under the map
$\pi$ from the preceding paragraph, these functions do \emph{not} correspond
to harmonic functions on~$\BB^d$ (or~$\bcv$).\footnote{And~there is little
reason why they should, because in dimension greater than 1 biholomorphic
maps do not preserve harmonicity in general (though they preserve holomorphy).}
In~particular, such~pushforwards to $\BB^d$ of harmonic functions on $\cB$ 
can be multi-valued at the origin.

With these prerequisites, we~now have the Poisson operator, Szeg\"o projection,
defining function, etc., on~$\cB$ (denoted by $\KK$, $\Pi$, $\rho$ and so~on;
for~the two-variable defining function one~can take the pullback 
$\rho((\CC z,cz),(\CC w,qw))=1-\spr{cz,qw}$ under $\pi$ of the two-variable
defining function $\rho(x,y)=1-\spr{x,y}$ for~$\BB^d$), as~well as the
analogous objects on~$\Omega=\cB\cap\cL_V$ (denoted $\KK_V$, $\Pi_V$,
$\rho_V=\rho|_{\cL_V\times\cL_V}$, and so~on), together with the restriction
operator $R_V$ from $\cB$ to~$\Omega$. (The~role of the ``ambient'' manifold
$W$ from Section~2 is played by $\cL$ for $\cB$ and by $\cL_V$ for~$\Omega$.)
Note~that from the above biholomorphism $\pi$ of 
$\cB\setminus\{\text{zero section}\}$ onto $\BB^d\setminus\{0\}$ sending
$\Omega\setminus\{\text{zero section}\}$ onto $\bcv\setminus\{0\}$,
it~follows that $\Omega$ (or,~more precisely,~$\cL_V$) intersects
$\partial\cB$ transversally (since $V$ intersects $\partial\BB^d$
transversally, thanks to the homogeneity of~$V$).

The~biholomorphism $\pi$ also identifies the spaces $\adac(\BB^d)$
with~$\adah(\cB)$, for~any $\alpha\in\RR$. Namely, using $\pi$
to transport the Lebesgue measure on $\CC^d\setminus\{0\}$ to
$\cL\setminus\{$zero section$\}$ (the~resulting volume element
on $\cL$ actually coincides with the one induced by the Hermitian
metric mentioned in the penultimate paragraph, provided the Fubini-Study
metric has been normalized so that $\CC P^{d-1}$ has volume~one),
the~fact that $\pi$ is a biholomorphism (except for the zero section
being sent to the origin, but these are both of measure~zero)
implies that $f\mapsto f\circ\pi$ is~a~unitary map from $L^2(\BB^d)$
onto $L^2(\cB)$ taking $L^2\hol(\BB^d)$ unitarily onto~$L^2\hol(\cB)$.
Similarly, since we are taking for the defining function $\rho$ on~$\cB$
the pullback under $\pi$ of the standard defining function $\rho_\BB(z)
=1-|z|^2$ on~$\BB^d$, the~last map acts unitarily from $L^2(\BB^d,
\rho_\BB^\alpha)$ onto $L^2(\cB,\rho^\alpha)$ for any $\alpha\in\RR$,
taking $\ada(\BB^d)$ unitarily onto $\ada(\cB)$ for any $\alpha>-1$.
By~the same argument with $\partial\BB^d$ and $\partial\cB$ in the places
of $\BB^d$ and~$\cB$, respectively, the~Hardy space $H^2(\partial\BB^d)$
is mapped by $\pi$ unitarily onto the Hardy space on~$\partial\cB$,
and the observation in the preceding sentence means that (cf.~(\tEE))
the \gto/ $T_{\Lambda^{\BB}_{\rho_{\BB}^\alpha}}$, $\Lambda^\BB_w=\KK^*_\BB
w\KK_\BB$, on~$\partial\BB^d$ is mapped by composition with~$\pi$ to the
\gto/ $T_{\Lambda_{\rho^\alpha}}$, $\Lambda_w=\KK^*w\KK$, on~$\partial\cB$
(even~though $\Lambda^\BB_{\rho_\BB^\alpha}$ is~\emph{not} in general mapped
into $\Lambda_{\rho^\alpha}$, because, as~was already remarked~above,
$\pi$~does not preserve harmonic functions). 
Similarly, repeating the arguments from the beginning of Section~4,
one~can check that the various norms on $\adah(\BB^d)$, $\alpha\in\RR$,
discussed there correspond under the composition with the biholomorphism
$\pi$ to norms of the form (\tQQ) with $T_Y$ an appropriate \gto/
on~$\partial\BB^d$. Quite generally, one~can arrive at this conclusion also
using the characterization of \gto/s as FIOs with the canonical relation
$\diag\Sigma$ that commute with~$\Pi$: indeed, as~$\pi$~is a biholomorphism,
the~composition with it is a FIO whose canonical relation takes
$\Sigma_{\partial\BB}$ isomorphically onto~$\Sigma_{\partial\cB}$, and,
as~we~have seen~above, intertwines the Szeg\"o projections on $\partial\BB^d$
and~$\partial\cB$; it~follows that the conjugation $T\mapsto\pi^{-1}T\pi$
is~a~symbol-preserving map from $I^m(\partial\BB^d,\partial\BB^d,\diag\Sigma
_{\partial\BB})$ onto $I^m(\partial\cB,\partial\cB,\diag\Sigma_{\partial\cB})$
that sends operators commuting with $\Pi_\BB$ into those commuting with~$\Pi$;
by~the above criterion, it~is thus an isomorphism from \gto/s
on~$\partial\BB^d$ onto those on~$\partial\cB$.

From~now~on, we~will therefore simply identify the spaces $\adac(\BB^d)$,
$\adah(\BB^d)$, and~so~on, discussed in the beginning of Section~4,
with the corresponding pullbacks $\adah(\cB)$ via $\pi$ on~$\cB$, 
knowing that the inner products in the latter are again of the form (\tQQ)
with $\partial\BB^d$ replaced by~$\partial\cB$ and $T_Y$ a suitable \gto/
on~$\partial\cB$; again, this includes in particular also the pullback 
$A^2_{-d,*}(\cB)$ under $\pi$ of the Drury-Arveson space $A^2_{-d}$ on~$\BB^d$.

\demo{Proof of Main Theorem} As~before, the~restriction of $R_V$ to
$\MMp=(\Ker R_V)^\perp\subset\adac(\BB^d)\equiv\adah(\cB)$ is a continuous
injective map of $\MMp$ into~$\ako$, which is now no longer onto in general
but by Corollary~{\PZ} its range is (closed~and) of~finite codimension.
Let~again $T_X$ be a positive selfadjoint elliptic \gto/ of order $-\alpha-k-1$
on $\pOm$ so that $\sprb{T_X\gamma_V f,\gamma_V g}$ is~the inner product
in~$\ako$, and $T_Y$ a positive selfadjoint elliptic \gto/ of order $-\alpha-1$
on $\partial\cB$ such that $\spr{T_Y\gamma f,\gamma g}_{\partial\cB}$ is~the
inner product in~$\adah(\cB)\equiv\adac(\BB^d)$. The~composed~map
$$ T = T_X^{1/2} \gamma_V R_V = T_X^{1/2} \rdv \gamma: \adac\to\hdva   $$
then satisfies $\Ker T=\MM$, maps boundedly $\adac$ into~$\hdva$,
and induces an isomorphism from $\MMp$ onto a (closed) subspace $\cN$ in
$\hdva$ of finite codimension. Let~$T^*:\hdva\to\adac$ be its adjoint.
By~abstract operator theory, the~operator $G:=(TT^*)|_{\cN} \oplus
I_{\cN^\perp}$ on $\cN\oplus\cN^\perp=\hdva$ is~then invertible
and for all $f,g\in\adac$
$$ \spr{P_\MMp f,g}_{\adac} = \sprb{G^{-1}Tf,Tg}.   \tag\tXM   $$
(Indeed, the~restriction $\tau$ of $T$ to $\MMp=(\Ker T)^\perp$ is injective
and maps onto $\Ran T=\cN$, hence is invertible as an operator from $\MMp$
onto~$\cN$; as~$TT^*=\tau\tau^*$ and $\cN=\Ran T=\Ran TT^*=(\Ker TT^*)^\perp$,
the~invertibility of $G$ follows. As~for~(\tXM), both sides vanish if $f$ or
$g$ belongs to~$\MM=\Ker T$; while for $f,g\in\MMp$, the~right-hand side
reduces to $\spr{(\tau\tau^*)^{-1}\tau f,\tau g}=\spr{f,g}$.)

Arguing as in the proof of Theorem~{\PY}, we~see from Proposition~{\PR} that
$TT^*$ is an elliptic \gto/ on $\pOm$ of order~0. Let~$T_H$ be its parametrix
(guaranteed by the property (P6) of \gto/s); thus $T_H$ is of order~0
(hence~bounded) and $TT^*T_H-I$ is a \gto/ of order $-\infty$.
By~Proposition~{\PC}, $TT^*T_H-I\in\SS^p$ for any $p>0$;
for~brevity, let~us temporarily denote an operator (not~necessarily the same
one at each occurrence) belonging to $\bigcap_{p>0}\SS^p$ by~$\cc$.
Since $G-TT^*=\cc$ (as~$I_{\cN^\perp}$ has finite rank), we~thus have
$GT_H=(TT^*+\cc)T_H=TT^*T_H+\cc=I+\cc$, whence $T_H=G^{-1}(I+\cc)=G^{-1}+\cc$.
Noting again that
$$ TM_{z_j} = T_X^{1/2}\gamma_V R_V M_{z_j} = T_X^{1/2}\gamma_V M_{z_j} R_V
 = T_X^{1/2} T_{z_j} \gamma_V R_V = T_X^{1/2}T_{z_j}T_X^{-1/2}T ,   $$
we~thus get from (\tXM) for any $f,g\in\MMp$,
$$ \align
\spr{S_j f,g}_\MMp &= \spr{P_\MMp M_{z_j} f,g}_{\adac}  \\
&= \sprb{G^{-1}T M_{z_j}f,Tg}  \\
&= \sprb{G^{-1}T_X^{1/2}T_{z_j}T_X^{-1/2}Tf,Tg},   \endalign   $$
that~is,
$$ S_j = \tau^* G^{-1} T_X^{1/2} T_{z_j} T_X^{-1/2} \tau
 = \tau^*(T_H+\cc) T_X^{1/2} T_{z_j} T_X^{-1/2} \tau
 = \tau^* \cT_j \tau + \cc   $$
(since $T_X^{1/2} T_{z_j} T_X^{-1/2}$, being a~\gto/ of order~0,
as~well as~$\tau$ are bounded), where
$$ \cT_j := T_H T_X^{1/2} T_{z_j} T_X^{-1/2}   $$
is~a~\gto/ of order~0. For any $j,l=1,\dots,d$, we~thus obtain analogously
as~before
$$ \align
[S_j,S^*_l] &= [\tau^* \cT_j \tau+\cc, \tau^* \cT^*_l \tau+\cc] 
 = [\tau^* \cT_j \tau, \tau^* \cT^*_l \tau] +\cc   \\
&= \tau^* (\cT_j[\tau\tau^*,\cT^*_l] + [\cT_j,\cT^*_l]\tau\tau^*
 + \cT^*_l[\cT_j,\tau\tau^*]) \tau + \cc.
  \endalign   $$
Once~again, the~last three commutators are \gto/s of order~$-1$ (or~less)
by~(\tCM), hence belong to $\SS^p$ for all $p>n$ by Proposition~{\PC},
and since the Schatten class $\SS^p$ forms an ideal in the algebra of
all bounded operators, we~get $[S_j,S^*_l]\in\SS^p$ $\forall p>n$,
proving the main theorem.   \qed   \enddemo

Again, we~have in fact proved that even $[S_j,S^*_l]\in\SS^{n,\infty}$.
Using the machinery of~\cite{\EZ}, it~is not difficult to give e.g.~a~formula
for the Dixmier trace of~$[S_j,S^*_l]^n$.

\medskip

We~finally conclude by observing that our results can be extended also to the
case when $V$ is a disjoint union of smooth submanifolds (away from the origin)
of possibly different dimensions.

\proclaim{Theorem~\prodef\LA} Let $V_1,\dots,V_m$ be homogeneous varieties
in $\CC^d$ such that $V_j\setminus\{0\}$ is a complex submanifold of $\CC^d
\setminus\{0\}$ of dimension~$n_j$, $j=1,\dots,m$, and $V_j\cap V_k=\{0\}$
for $j\neq k$. Let $\alpha\in\RR$ and let $\MM$ be the subspace in~$\adac$
$($or~in $\ada$ if~$\alpha>-d-1)$ of functions that vanish on~$V\cap\BB^d:=
\bigcup_{j=1}^m V_j\cap\BB^d$. Then $[S_j,S^*_k]\in\SS^q$, $j,k=1,\dots,d$,
for all $q>\max(n_1,\dots,n_m)$.    \endproclaim

\demo{Proof} Let $R_{V_j}$ be the restriction operator for $\cL_{V_j}=
\pi^{-1}(V_j)$ and $T_{X_j}$ be the positive selfadjoint elliptic \gto/
of order $-\alpha-(d-n_j)-1$ on $\pOm_j$, $\Omega_j:=\cL_{V_j}\cap\cB$,
as~in the preceding proof for $V_j$ in the place of the $V$ there,
$j=1,\dots,m$; and let also $T_Y$ be as in the preceding proof.
Denote by $T$ the the column block matrix with entries $T_j:=T_{X_j}^{1/2}
\gamma_{V_j}R_{V_j}=T_{X_j}^{1/2}R_{\partial V_j}\gamma$, $j=1,\dots,m$;
thus $T$ acts continuously from $\adac(\BB^d)\equiv\adah(\cB)$ into the
Hilbert space direct sum $\cH:=\bigoplus_{j=1}^m H^2(\pOm_j)$, and~as before
$\Ker T=\MM$. We~have seen in the preceding proof that each $T_jT^*_j$
is a Fredholm operator on~$H^2(\pOm_j)$, by~Proposition~{\PR};
on~the~other~hand, from the proof of that proposition we~also see that 
$T_jT^*_k$ for $j\neq k$ is an operator from $H^2(\pOm_k)$ into $H^2(\pOm_j)$
whose Schwartz kernel is in $C^\infty(\pOm_j\times\pOm_k)$: namely, the~latter
kernel is the restriction of $(T_{X_j}^{1/2}\otimes\widetilde T_{X_k}^{1/2})
T_Y^{-1}S_x(y)$ to $x\in\pOm_j$ and $y\in\pOm_k$ (where the tensor product
notation means that $T_{X_j}^{1/2}$ acts on the $x$ variable and $\widetilde
T_{X_k}^{1/2}$ on the $y$~variable, and for any operator $A$ one defines
$\widetilde A f:=\overline{A^*\overline f}$, with bar denoting complex
conjugation), and $T_Y^{-1}S_x(y)$ has singularities
only on the diagonal $x=y$ by~(\tSZ) while $\pOm_j\cap\pOm_k=\emptyset$
by~hypothesis (also $T_{X_j}^{1/2}\otimes\widetilde T_{X_k}^{1/2}$ maps
$C^\infty(\pOm_j\times\pOm_k)$ into itself, in~view of the way \gto/s act
on~Sobolev spaces). 
Thus $T_jT^*_k$ is a smoothing operator for $j\neq k$, and hence belongs
to all~$\SS^p$, $p>0$. Denoting temporarily by $D$ the $m\times m$ block
matrix with $T_jT^*_j$, $j=1,\dots,m$, on~the main diagonal and zeroes
elsewhere, we~therefore~have
$$ TT^* = D + \cc   \tag\tDC   $$
where $\cc$ has the same meaning as in the preceding proof. It~follows that
$TT^*$ is~again a Fredholm operator, and the restriction $\tau$ of $T$ to
$\MMp$ is an isomorphism of $\MMp$ onto the (closed) subspace $\cN:=\Ran T$
in $\cH$ of finite codimension; the~operator $G:=(TT^*)|_{\cN}\oplus
I_{\cN^\perp}$ on $\cN\oplus\cN^\perp=\cH$ is invertible and for all
$f,g\in\adac$,
$$ \spr{P_{\MMp}f,g}_{\adac} = \spr{G^{-1}Tf,Tg}_{\cH}.   $$
Let~$T_{H_j}$ be a parametrix for $T_jT^*_j$, $j=1,\dots,m$ (this is an
elliptic \gto/ of order 0 on~$\pOm_j$), and $H$ the $m\times m$ block
matrix with $T_{H_j}$ on the main diagonal and zeroes elsewhere.
By~(\tDC), $TT^*H-I\in\cc$ and arguing as before, we~get $H=G^{-1}+\cc$ and
$$ S_l = \tau^* \cT_l \tau + \cc ,   $$
where $\cT_l$ is the $m\times m$ block matrix with $T_{H_j}T_{X_j}^{1/2}
T_{z_l}T_{X_j}^{-1/2}$, $j=1,\dots,m$, on~the main diagonal and zeroes
everywhere else; thus $\cT_l$ is a direct sum over $j$ of \gto/s of order~0
on~$H^2(\pOm_j)$. For~any $k,l=1,\dots,d$ we thus again~get
$$ \align
[S_j,S^*_l] &= [\tau^* \cT_j \tau+\cc, \tau^* \cT^*_l \tau+\cc] 
 = [\tau^* \cT_j \tau, \tau^* \cT^*_l \tau] +\cc   \\
&= \tau^* (\cT_j[\tau\tau^*,\cT^*_l] + [\cT_j,\cT^*_l]\tau\tau^*
 + \cT^*_l[\cT_j,\tau\tau^*]) \tau + \cc  \\
&= \tau^* (\cT_j[D,\cT^*_l] + [\cT_j,\cT^*_l] D
 + \cT^*_l[\cT_j,D]) \tau + \cc .   \endalign   $$
The~last three commutators are $m\times m$ block matrices with \gto/s of
order~$-1$ on the main diagonal and zeroes elsewhere, hence belong to
$\bigoplus_{j=1}^m \SS^{p_j}(H^2(\pOm_j))$ for all $p_j>n_j$, $j=1,\dots,m$,
by~Proposition~{\PC}, so~to $\SS^p(\cH)$ if $p>n:=\max(n_1,\dots,n_m)$.
Since the Schatten classes form an ideal and $T$ is bounded, 
the~theorem follows. (And~again, in~fact $[S_k,S^*_l]\in\SS^{n,\infty}$.)
\qed   \enddemo 

We~remark that once we know that the block matrix operator~$TT^*$ is~Fredholm,
and, hence, $T$~has finite-codimensional range, it~follows by abstract operator
theory that both $T$ and $T^*$ have closed ranges; but the range of the row
block operator $T^*$ is precisely the sum $\MM_1^\perp+\dots+\MM_m^\perp$,
where $\MM_j$ is the subspace in $\adac$ of functions vanishing on~$V_j$.
By~Proposition~3.4(2) in Kennedy and Shalit~\cite{\KS}, the~closedness
of this sum together with our Main Theorem imply that $[S_k,S^*_l]\in\SS^p$,
$k,l=1,\dots,d$, if~$p>\max(n_1,\dots,n_m)$, thus yielding an
alternative way to conclude the last proof.

An~elementary argument using the decomposition of homogeneous varieties into
their irreducible components shows that every homogeneous variety in $\CC^d$
which is smooth at each of its non-zero points is a finite union of varieties
$V_j$ as in Theorem~\LA. Therefore Theorem~{\LA} generalizes our main theorem
to the case of arbitrary homogeneous varieties $V$ in $\CC^d$ that are smooth
outside the origin. That~is, the~refinement of the Arveson Conjecture as
formulated by Douglas holds in this case. (Note that the dimension of the
analytic set $V = \bigcup_{j=1}^m V_j$ at the origin is given by $\dim_0 V
=\max(n_1,...,n_m)$, see Section~5.3 in~\cite{\GrauRem}.)

\Refs
\widestnumber\key{99}
\konectrue

\refd\ArvP {W. Arveson: {\it $p$-summable commutators in dimension~$d$,\/}
 J.~Operator Theory {\bf 54} (2005), 101--117}

\refd\BdMA {L. Boutet de Monvel: {\it Boundary problems for pseudo-differential
operators,\/} Acta Math. {\bf 126} (1971), 11--51}

\refd\BdMi {L. Boutet de Monvel: {\it On the index of Toeplitz operators 
in several complex variables,\/} Invent. Math. {\bf 50} (1979), 249--272}

\refd\BdMG {L. Boutet de Monvel, V. Guillemin: {\it The spectral theory of
Toeplitz operators,\/} Ann. Math. Studies, vol.~99, Princeton University 
Press, Princeton, 1981}

\refd\BdMS {L. Boutet de Monvel, J. Sj\"ostrand: {\it Sur la singularit\'e des
noyaux de Bergman et de Szeg\"o,\/} Ast\'erisque {\bf 34--35} (1976), 123--164}
 
\refd\BdMT {L. Boutet de Monvel: {\it Symplectic cones and Toeplitz
operators,\/} Multidimensional Complex Analysis and Partial Differential
Equations, pp.~15--24, Contemporary Math., vol.~205, Amer. Math. Soc.,
Providence, 1997}

\refd\Bea {F. Beatrous: {\it Estimates for derivatives of holomorphic
functions in pseudoconvex domains,\/} Math. Z. {\bf 191} (1986), 91--116}

\refd\Doug {R. Douglas: {\it A~new kind of index theorem,\/} Analysis,
geometry and topology of elliptic operators, pp.~369-382, World Sci. Publ.,
Hackensack, 2006}

\refd\DW {R.G. Douglas, K.Wang: {\it Essential normality of cyclic submodule
generated by any polynomial,\/} preprint, 2011, arXiv:1101.0774}

\refd\Ejfa {M. Engli\v s: {\it Toeplitz operators and weighted Bergman
kernels,\/} J.~Funct. Anal. {\bf 255} (2008), 1419--1457}

\refd\Epay {M. Engli\v s: {\it Weighted Bergman kernels for logarithmic
weights,\/} Pure Appl. Math. Quarterly (Kohn special issue) {\bf6} (2010),
781--813} 

\refd\E {M. Engli\v s: {\it Analytic continuation of weighted Bergman
kernels,\/} J.~Math. Pures Appl. {\bf 94} (2010), 622--650}

\refd\EZ {M. Engli\v s, G. Zhang: {\it Hankel operators and the Dixmier
trace on strictly pseudoconvex domains,\/} Docum. Math. {\bf 15} (2010),
601--622}

\refd\FX {Q. Fang, J. Xia: {\it Essential normality of polynomial-generated
submodules: Hardy space and beyond,\/} preprint, 2011, available~at 
http://www.acsu.buffalo.edu/$^\sim$jxia/Preprints/poly.pdf}

\refd\GrauRem {H. Grauert, R.Remmert, {\it Coherent analytic sheaves,\/}
Springer, Berlin, 1984}

\refd\GrSj {A. Grigis, J. Sj\"ostrand, {\it Microlocal analysis for
diferential operators,\/} London Math. Soc. Lecture Notes, vol.~196,
Cambridge Univ. Press, Cambridge, 1994}

\refd\GW {K. Guo and K. Wang, {\it Essentially normal Hilbert modules and
K-homology,\/} Math. Ann. {\bf 340} (2008), 907--934}

\refd\Hcpam {L. H\"ormander: {\it Pseudo-differential operators,\/} Comm.
Pure Appl. Math. {\bf18} (1965), 501--517}

\refd\Horm {L. H\"ormander, {\it The analysis of linear partial differential
operators, vol.~I--IV,\/} Grund\-leh\-ren der mathematischen Wissenschaften,
Springer-Verlag, 1985}

\refd\KS {M. Kennedy, O.M. Shalit: {\it Essential normality and the
decomposability of algebraic varieties,\/} New~York J. Math. {\bf 18}
(2012), 877--890}

\refd\LM {J.-L. Lions, E. Magenes, {\it Probl\`emes aux limites non
homog\`enes et applications,\/} vol.~1, Dunod, Paris, 1968}

\refd\MeSj {A. Melin, J. Sj\"ostrand: {\it Fourier integral operators
with complex valued phase functions,\/} Fourier integral operators and
partial differential equations, Lecture Notes Math., vol.~459, pp.~120-223,
Springer Verlag, Berlin-Heidelberg, 1975}

\refd\Sh {O.M. Shalit: {\it Operator theory and function theory in
Drury-Arveson space and its quotients,\/} preprint, 2013, arxiv:1308.1081}

\refd\Strich {R. Strichartz: {\it A~functional calculus for elliptic
pseudo-differential operators,\/} Amer. J. Math. {\bf 94} (1972), 711--722}

\refd\Trev {F. Treves, {\it Introduction to pseudodifferential and Fourier
integral operators,\/} Plenum, New~York, 1980}

\endRefs

\enddocument

\bye